\input amstex
\documentstyle{amsppt}
\document
\magnification=1200
\NoBlackBoxes
\nologo
\pageheight{18cm}

\bigskip

\centerline{\bf INVERTIBLE COHOMOLOGICAL FIELD THEORIES}

\smallskip

\centerline{\bf AND WEIL--PETERSSON VOLUMES}

\bigskip

\centerline{Yuri I. Manin$^{1)}$ and Peter Zograf$^{1),2)}$}

\medskip

\centerline{\it 1) Max-Planck-Institut f\"ur Mathematik, Bonn, Germany}

\centerline{\it 2) Steklov Mathematical Institute, St.Petersburg, Russia}

\bigskip

{\bf Abstract.} We show that the generating function for the
higher Weil--Petersson volumes of the moduli spaces of stable
curves with marked points can be obtained from  Witten's
free energy by a change of variables given by Schur polynomials.
Since this generating function has a natural extension
to the moduli space of invertible Cohomological Field Theories,
this suggests the existence of a ``very large phase space'',
correlation functions on which include Hodge integrals
studied by C.~Faber and R.~Pandharipande. From this formula we derive
an asymptotical expression for the Weil--Petersson volume
as conjectured by C.~Itzykson. 
We also discuss a topological interpretation of the genus expansion
formula of Itzykson--Zuber, as well as a related bialgebra acting 
upon quantum cohomology as a complex version of the classical path groupoid.

\bigskip

{\bf 0. Introduction and summary.} The aim of this paper is to record some
progress in understanding intersection numbers on moduli
spaces of stable pointed curves and their
generating functions. Continuing the study started in [KoM],
[KoMK] and pursued further in [KaMZ], [KabKi], we work
with Cohomological Field Theories (CohFT), which is the same as
cyclic algebras over the modular operad $H_*M_{g,n}=H_*(\overline{M}_{g,n+1},K)$
in the sense of [GK] (here $K$ is a field of characteristic zero).
Such algebras form a category with symmetric tensor product and identical
object. Isomorphism classes of (slightly rigidified) invertible objects 
form $K$--points of an infinite--dimensional
abelian algebraic group which can be called {\it the Picard
group} of the respective category. This group can be identified 
with the product of
$K^*$ and a vector space $L$ consisting of certain families
of elements in all $H^*(\overline{M}_{g,n},K)^{S_n}.$ The space $L$ is
naturally graded by (co)dimension. 

\smallskip

(A) The first
problem is the calculation of $L$. To our knowledge, it is unsolved,
and we want to stress its importance. We know two independent elements
$\kappa_a$ and $\mu_a$ in each odd codimension $a\ge 1$, and one
element  $\kappa_a$ in each even codimension $a\ge 2$. There are
only two divisorial classes in $L$, namely, $\kappa_1$ and $\mu_1$.
As it follows from the results of [KaMZ], the genus zero restriction of
$L$  is generated by the $\kappa$-classes, and in genus one there is just
one additional class $\mu_1$ (see [KabKi]). 

\smallskip

(B) The second problem
is a description of a formal function on $L$ (or $K^*\times L)$,
the potential of the respective invertible CohFT.

\smallskip

It turns out that evaluating this function {\it on the subspace
generated by $\kappa_a$ only}, we get essentially E.~Witten's
generating function for $\tau$--intersection numbers
(total free energy of two dimensional gravity)
but written in a different coordinate system: the
respective nonlinear coordinate change is given by Schur
polynomials. The proof uses the
explicit formula for higher Weil--Petersson volumes of arbitrary genus
derived in [KaMZ]. This completes the
calculation of the total generating function for these
volumes: the genus zero
component was calculated in [KaMZ], and the results
of [IZu] about the genus expansion of the total free energy 
now allow us to express all higher genus contributions through the genus
zero one. This is a version of the argument
given in [Zo], see also Theorem 4.1 below. It also leads to a complete
proof of Itzykson's conjecture about asymptotics of
the classical Weil--Petersson volumes of $\overline{M}_{g,n}$,
see [Zo], formula (7):
$$
\frac{\langle\kappa_1^{3g-3+n}\rangle}{n! (3g-3+n)!}\, =\, 
C^n\, n^{-1+\frac{5}{2}(g-1)}(B_g+O(1/n)),\quad n\to\infty 
$$
For genus zero it was proved in [KaMZ]. The constant $C$ does not 
depend on $g$. It was mentioned  in [KaMZ] that it is expressible
via the first zero of the Bessel function $J_0$.
We supply here a detailed proof and extend it to all genera.

\smallskip

Looking at the  identification of generating functions
from the other side, we realize that the complete potential 
as a function on the moduli space
of invertible CohFT's is a natural infinite dimensional extension of
the total free energy. It would be interesting to derive
differential equations for this function
extending the Virasoro constraints and higher KdV equations.

\smallskip

Finally, we discuss several related constructions
and analogies: 

\smallskip

(C) The cohomology space of any projective algebraic smooth
manifold $V$ carries a canonical structure of CohFT (quantum
cohomology).
The identity object of the category CohFT is the quantum cohomology
of a point. It is tempting to interpret any other
invertible CohFT, $T$, as  generalized quantum cohomology
of a point, and to define {\it the generalized
quantum cohomology of $V$} as the cyclic modular algebra
$H^*(V)\otimes T.$ The potential of this theory
restricted to the $\kappa$--subspace can be derived from
the potential of $V$ with gravitational descendants.
Of course, it is an {\it ad hoc} prescription which
must be replaced by a more geometric construction. 

\smallskip

This sheds new light on the mathematical nature of the {\it large
phase space} (see [W]): according to the
comments (B), it must be a part of a {\it very
large phase space} which is a suitably
interpreted moduli space of generalized quantum cohomology.

\smallskip

(D) The structure of the operad $H^*(\overline{M}_{g,n},K)$
restricted to the $n=2$ part of it gives rise to an interesting bialgebra.
A part of its primitive elements can be obtained from $L.$
It would be interesting to describe all of them.
The action of this bialgebra on $H^*(V)$ can be derived from
a complexified version of the path groupoid of $V.$
This raises several problems of which the most interesting is
probably the construction of the quantum motivic fundamental
group, upon which this groupoid might act.

\medskip

{\bf 1. The space $L.$} By definition, $L\subset \prod_{g,n}
H^*(\overline{M}_{g,n},K)$ is the linear space formed
by all families $l=\{l_{g,n}\in H^*(\overline{M}_{g,n},K)\}$
satisfying the following coditions:

\smallskip

(i) $l_{g,n}\in H^{ev}(\overline{M}_{g,n},K)^{S_n}.$

\smallskip

(ii) For every boundary morphism 
$$
b:\, \overline{M}_{g_1,n_1+1}\times\overline{M}_{g_2,n_2+1}\to
\overline{M}_{g,n} 
\eqno(1)
$$ 
with $g=g_1+g_2, n_1+n_2=n,$ we have
$$
b^*(l_{g,n})=l_{g_1,n_1+1}\otimes 1+1\otimes l_{g_2,n_2+1}.
\eqno(2)
$$  

\smallskip

(iii) For any boundary morphism
$$
b^{\prime}:\ \overline{M}_{g-1,n+2}\to \overline{M}_{g,n}
\eqno(3)
$$
we have
$$
b^{\prime *}(l_{g,n}) = l_{g-1,n+2}.
\eqno(4)
$$

\smallskip

Recall that the boundary morphisms $b$ (resp. $b^{\prime}$)
glue together
a pair of marked points situated on different
connected components (resp. on the same one) of the 
appropriate universal curves.

\medskip

{\bf 1.1. Proposition.} {\it (a) $L$ is graded by codimension in
the following sense: 
put $L^{(a)}\,=L\,\cap\,\prod_{g,n}
H^{2a}(\overline{M}_{g,n},K)$, then $L=\prod_{a\ge 1}L^{(a)}.$

\smallskip

(b) Let $(H,h)$ be a one--dimensional vector space
over $K$ endowed with a metric $h$ and an even basic vector
$\Delta_0$ such that $h(\Delta_0,\Delta_0)=1.$
For any structure of the invertible CohFT on $(H,h)$
given by the maps $I_{g,n}:\,H^{\otimes n}\to H^*(\overline{M}_{g,n},K)$
put $c_{g,n}=I_{g,n}(\Delta_0^{\otimes n}).$ 
Let $P$ be the group of the isomorphism classes of
such theories with respect to the tensor multiplication
(isomorphisms should identify the basic vectors). Then
the following map $K^*\times L\to P$ is a group
isomorphism:
$$
(t,\{l_{g,n}\})\mapsto \{c_{g,n}=t^{2g-2+n}\,\roman{exp}\,(l_{g,n})\}.
\eqno(5)
$$  
}

\smallskip

The proof is straightforward: cf. [KoMK], sec. 3.1, where most of it
is explained for genus zero. Notice that
the statement on the grading means that if one puts 
$l_{g,n}=\sum_a l^{(a)}_{g,n}$ with $l^{(a)}_{g,n}\in 
H^{2a}(\overline{M}_{g,n},K),$ then 
$l\in L$ is equivalent to $l^{(a)}\in L$ for all $a,$
where $l^{(a)}=\{l^{(a)}_{g,n}\}.$ 

\medskip

{\bf 2. Classes $\kappa_a$ and $\mu_a$.} Put
$$
\kappa_a=\{\,\kappa_{g,n;\,a}=\pi_*(c_1(\omega_{g,n}(\sum_{i=1}^n x_i
))^{a+1}) \in H^{2a}(\overline{M}_{g,n},K)^{S_n}\}
\eqno(6)
$$
where $\omega_{g,n}$ is the relative dualizing sheaf of the
universal curve $\pi =\pi_{g,n}:\,C_{g,n}\to\overline{M}_{g,n},$
and $x_i$ are the structure sections of $\pi .$

\smallskip

Furthermore, with the same notation put
$$
\mu_a=\{\,\mu_{g,n;\,a}=\roman{ch}_a (\pi_* \omega_{g,n})\,\}.
\eqno(7)
$$
Classes $\mu_{g,n;a}$ vanish for all even $a$: see [Mu].

\smallskip

Both $\kappa_a$ and $\mu_a$ belong to $L^{(a)}$: for $\kappa$
this was noticed, e.g., in [AC], and for $\mu$ in [Mu].
These elements are linearly independent for odd $a$ 
because, for example, $\mu_{g,n;\,a}$ is lifted from $\overline{M}_{g,0}$
for $g\ge 2$, or from $\overline{M}_{1,1},$ whereas
$\kappa_{g,n;\,a}$ is not.

\medskip

{\bf 2.1. Remark.} The classes $\kappa_{g,n;\,a}$ are denoted
$\omega_{g,n}(a)$ in [KaMZ]. We use this opportunity to clarify
one apparent notational ambiguity in that paper. Namely,
the initial definition of $\omega_{g,n}(a)$ given by (0.1) in [KaMZ]
coincides with our (6). On the other hand, specializing
the formula (2.4) of [KaMZ] for $\omega_{g,n}(a_1,\dots ,a_p)$ 
to the case $p=1$ we obtain a formally different expression
$$
\omega_{g,n}(a)=\pi_*(c_1(x_{n+1}^*\omega_{g,n+1})^{a+1}),
$$
where $C_{g,n}$ is canonically identified with $\overline{M}_{g,n+1}$
and $x_{n+1}\!:\, \overline{M}_{g,n+1}\rightarrow C_{g,n+1}$ is the respective section.
Actually, the two definitions coincide because
there exists a canonical isomorphism on $C_{g,n}$
$$
\omega_{g,n}(\sum_{i=1}^n x_i)\simeq x_{n+1}^*\omega_{g,n+1}.
$$

\medskip

{\bf 3. Three generating functions.} Witten's {\it total
free energy} of two dimensional gravity
is the formal series
$$
F(t_0,t_1,\dots )=\sum_{g=0}^{\infty} F_g(t_0,t_1,\dots )=
\sum_{g=0}^{\infty}\quad\sum_{\sum (i-1)l_i=3g-3}\langle\tau_0^{l_0}\tau_1^{l_1}\dots
\rangle\,\prod_{i=0}^{\infty}\frac{t_i^{l_i}}{l_i!}.
\eqno(8)
$$
Here $\langle\tau_0^{l_0}\tau_1^{l_1}\dots\rangle$
is the intersection number defined as follows. Consider
a partition $d_1+\dots + d_n=3g-3+n$ such that
$l_0$ of the summands $d_i$ are equal to 0, $l_1$ are equal to 1, etc, and put
$$
\langle\tau_0^{l_0}\tau_1^{l_1}\dots\rangle =
\langle\tau_{d_1}\dots \tau_{d_n}\rangle =
\int_{\overline{M}_{g,n}}\psi_1^{d_1}\dots\psi_n^{d_n}
\eqno(9)
$$
where $\psi_i= \psi_{g,n;\,i}= c_1(x_i^*\omega_{g,n})$
and $n=\sum l_i$.

\medskip

The generating function for {\it higher Weil--Petersson volumes} was 
introduced (in a slightly different form) in [KaMZ]:
$$
K(x,s_1, \dots )=\sum_{g=0}^{\infty} K_g(x,s_1, \dots )=
$$
$$
\sum_{g=0}^{\infty}\,\sum_{n=0}^{\infty}\; \sum_{|\bold{m}|=3g-3+n}
\langle \kappa_1^{m_1} \kappa_2^{m_2} \dots\rangle\,
\frac{x^n}{n!}\,\prod_{a=1}^{\infty}\frac{s_a^{m_a}}{m_a!}\, .
\eqno(10)
$$
Here and below $\bold{m}=(m_1,m_2,\dots m_{3g-3+n})$ with $|\bold{m}|=\sum_{a=1}^{3g-3+m} am_a$, and
$$
\langle \kappa_1^{m_1} \kappa_2^{m_2} \dots \rangle
= \int_{\overline{M}_{g,n}}
\kappa_{g,n;\,1}^{m_1}\kappa_{g,n;\,2}^{m_2}\dots 
\kappa_{g,n;\,3g-3+n}^{m_{3g-3+n}}
\eqno(11)
$$
(all integrals over the moduli spaces of unstable curves 
are assumed to be zero).

\medskip

Actually, (10) is a specialization of the general
notion of the the {\it potential} of a CohFT,
$(H,h;\{I_{g,n}\})$, which is
a formal series in coordinates of $H$ and
is defined as $\langle \roman{exp}\,(\gamma )\rangle$ where
$\gamma$ is the generic even element of $H$ and
the functional $\langle *\rangle$ as above is the
integration over fundamental classes. Considering variable
invertible CohFT's, we will get the potential as
a formal function on $H\times K^* \times L$ which is a
series in $(x, t, s_a, r_b, \dots ).$ Here $x$
stands for the coordinate on $H$ dual to $\Delta_0$, 
$t$ for the coordinate on $K^*$ as in (5), $s_a$ (resp. $r_b$)
are coordinates on $L$ dual to $\kappa_a$ (resp. $\mu_a$),
and dots stand for the unexplored part of $L$ (if there is any).
More precisely, potential of the individual theory (5) is
$$
\Phi (x, t) = \sum_g \Phi_g(x, t) =\sum_g \langle \roman{exp}\,(x\Delta_0)\rangle_g =
$$
$$
\sum_{g,n}\frac{\langle I_{g,n}(x^n\Delta_0^{\otimes n})\rangle}{n!}=
\sum_{g,n} \frac{x^n}{n!}\,\int_{\overline{M}_{g,n}}t^{2g-2+n}
\,\roman{exp}\,(l_{g,n}) \, .
$$
Now make $l_{g,n}$ generic, that is, put 
$$
l_{g,n}=\sum_{a,b,\dots}(s_a\kappa_{g,n;\,a}+r_b\mu_{g,n;\,b}
+\dots )
$$ 
Collecting the terms of the right dimension,
we finally get the series
$$
\Phi (x, t, s_1,s_2, \dots , r_1,r_2,\dots ) = 
\sum_g \Phi_g(x, t, s_1,s_2,\dots , r_1,r_2,\dots )=
$$
$$
\sum_{g,n} \frac{x^n}{n!}\,t^{2g-2+n}\sum
\langle \kappa_1^{m_1} \kappa_2^{m_2} \dots
\mu_1^{p_1} \mu_2^{p_2} \dots\rangle
\,\prod_{a=1}^{\infty}\frac{s_a^{m_a}}{m_a!}\,
\,\prod_{b=1}^{\infty}\frac{r_b^{p_b}}{p_b!}\,\dots
\eqno(12)
$$
where the inner summation is taken over $\bold{m},\bold{p}, \dots$ with
$|\bold{m}|+|\bold{p}|+\dots =3g-3+n.$

\smallskip

Putting here $t=1, r_b=\dots =0,$ we obtain (10).
The following Theorem 4.1 shows, that making an
invertible change of variables in (8) restricted
to $t_1=0$ we again can get (10). Thus (12)
can be considered as a natural 
infinite dimensional extension of (8).

\smallskip

Some information about the coefficients of (12) with
non--vanishing $p_b$ is obtained in [FP]. If
$L$ is not generated by the $\kappa$-- and $\mu$--classes,
these coefficients provide the natural generalization
of Hodge integrals of [FP].

\medskip

{\bf 4. Relation between generating functions.} The main result
of this section is:

\medskip

{\bf 4.1. Theorem.} {\it For every $g=0,1, \dots$ we have
$$
K_g(x,s_1,\dots )=F_g(t_0,t_1,\dots )|_{t_0=x,t_1=0,t_{j+1}=p_j(s_1,\dots ,s_j)}
$$
Here $p_j$ are  the Schur polynomials defined by
$$
1-\roman{exp}\,\left(-\sum_{i=1}^{\infty}\lambda^is_i\right)=
\sum_{j=1}^{\infty}\lambda^j p_j(s_1,\dots ,s_j).
$$
}

\medskip

{\bf Proof.} 
The statement of the theorem is a convenient reformulation of 
formula (2.17) of [KaMZ]:
$$
\frac{\langle \kappa_1^{m_1} \dots \kappa_{3g-3+n}^{m_{3g-3+n}}\rangle}
{m_1!\dots m_{3g-3+n}!}=
$$
$$
\sum_{k=1}^{\|\bold{m}\|}\frac{(-1)^{\|\bold{m}\|-k}}{k!}
\sum_{\bold{m}=\bold{m}^{(1)}+\dots +\bold{m}^{(k)}}
\frac{\langle\tau_0^n\tau_{|\bold{m}^{(1)}|+1}\dots \tau_{|\bold{m}^{(k)}|+1}\rangle}
{\bold{m}^{(1)}!\dots \bold{m}^{(k)}!} \, ,
\eqno(13)
$$
where $\|\bold{m}\|=\sum m_i,\;\;\bold{m}^{(i)}=(m_1^{(i)},\dots , m_{3g-3+n}^{(i)})\ne 0$,
$\bold{m}^{(i)}!=m_1^{(i)}!\dots m_{3g-3+n}^{(i)}!$\,.
To show this, fix $\bold{m}=(m_1, \dots , m_{3g-3+n})$ and compare the coefficients
at the monomial $x^n s_1^{m_1}\dots s_{3g-3+n}^{m_{3g-3+n}}$ 
in the series $K_g$ and $F_g(x,0,p_1,p_2,\dots)$.  
The contribution from the term 
$$
\langle\tau_0^{l_0}\tau_1^{l_1}\dots\rangle\,\prod \frac{t_j^{l_j}}{l_j!}
\eqno(14)
$$
to this coefficient
is nontrivial if and only if $l_0=n,\; l_1=0$ and there exists a partition 
$\bold{m}=\bold{m}^{(1)}+\dots +\bold{m}^{(k)}$ with $k=l_2+l_3+\dots$ such that for any $j=2,3,\dots$
there are exactly $l_j$ sequences $\bold{m}^{(i)}=(m^{(i)}_1,m^{(i)}_2,\dots)$ 
with $|\bold{m}^{(i)}|=j-1$. Then, obviously,
$\langle\tau_0^{l_0}\tau_1^{l_1}\dots\rangle =
\langle\tau_0^n\tau_{|\bold{m}^{(1)}|+1}\dots \tau_{|\bold{m}^{(k)}|+1}\rangle\,.$
Let us subdivide the set of $k=l_2+l_3+\dots$ sequences $\bold{m}^{(i)}$ into subsets
consisting of equal sequences and denote the cardinalities of these subsets
by $l_2^{(1)},l_2^{(2)},\dots, l_3^{(1)},l_3^{(2)},\dots, l_j^{(1)},l_j^{(2)},\dots,$
so that $\sum_a l_j^{(a)}=l_j$. 
Since the polynomials $p_j$ are explicitly given by
$$
p_j(s_1,\dots ,s_j)=-\sum_{|\bold{m}|=j}\prod_{i=1}^j (-1)^{m_i}\frac{s_i^{m_i}}{m_i!},
$$
we can rewrite the contribution from (14) corresponding to the partition
$\bold{m}=\bold{m}^{(1)}+\dots +\bold{m}^{(k)}$ as
$$
\langle\tau_0^n\tau_{|\bold{m}^{(1)}|+1}\dots \tau_{|\bold{m}^{(k)}|+1}\rangle\,\frac{x^n}{\prod l_j!}\,
\,\prod_{j\ge 2}\frac{l_j!}{\prod_a l_j^{(a)}!}\prod_{i=1}^k \left(-\prod_{b=1}^{3g-3+n}
(-1)^{m_b^{(i)}}\frac{s_b^{m_b^{(i)}}}{m_b^{(i)}!}\right)=
$$
$$
(-1)^{\|\bold{m}\|-k}\,\frac{1}{\prod_{j,a}l_j^{(a)}!}\,\, 
\frac{\langle\tau_0^n\tau_{|\bold{m}^{(1)}|+1}\dots \tau_{|\bold{m}^{(k)}|+1}\rangle}{\bold{m}^{(1)}!\dots \bold{m}^{(k)}!}
\,\,\frac{x^n s_1^{m_1}\dots s_{3g-3+n}^{m_{3g-3+n}}}{n!}.
\eqno(15)
$$
Since there are exactly $k!/\prod_{j,a}l_j^{(a)}!$ partitions
of $\bold{m}=(m_1, \dots , m_{3g-3+n})$ that are equivalent 
to the above partition  $\bold{m}=\bold{m}^{(1)}+\dots +\bold{m}^{(k)}$
under permutation of indices $1,\dots,k$, 
the coefficient in (15) coincides with
the corresponding term in (13), which completes the proof.

\medskip

{\bf 4.2. Corollaries.} It follows from Theorem 4.1 that
all known properties of the function $F$ can be automatically
translated to $K$. First, $F$ satisfies the Virasoro constraints
(or, equvalently, $\partial^2 F/\partial t_0^2$ is a solution to
the KdV hierarchy). The same is true for $K$
up to a change of variables.

\smallskip

Second, we can express each $K_g$ with $g\ge 1$ in terms of $K_0$
in the same way as it was done in [IZu] for $F.$ Let us remind the genus
expansion of $F$ obtained in [IZu] (cf. also [EYY]). Put
$u_0=F_0^{\prime\prime}$ where a prime denotes the derivative with respect to
$t_0$, and define a sequence of formal series $I_1, I_2, \dots$ by
$$
I_1=1-\frac{1}{u_0^{\prime}},\quad I_{k+1}= \frac{I_k^{\prime}}{u_0^{\prime}},
\quad k=1, 2, \dots
$$
Then $F_1=\frac{1}{24}\,\roman{log}\,u_0^{\prime}$, and for $g\ge 2$
$$
F_g=\sum_{\sum(i-1)m_i=3g-3} \langle \tau_2^{m_2}\tau_3^{m_3}\dots
\tau_{3g-2}^{m_{3g-2}}\rangle\,
u_0^{\prime\, 2g-2+\|\bold{m}\|}\,\prod_{k=2}^{3g-2}\frac{I_k^{m_k}}{m_k!}\, .
\eqno(16)
$$
Substituting here $t_0=x, t_1=0$ and $t_{i+1}=p_i(s_1,\dots ,s_i), i\ge 1$
we get the identical formula for $K_g.$

\medskip

{\bf 5. Topological interpretation of the genus expansion.} 
Differentiating $n$ times the both sides of (16) with respect to $t_0$
and using the definition of $I_k$, one easily gets that
$$
\frac{\partial^n F_g}{\partial t_0^n}=\sum_{|{\bold{m}}|=3g-3+n}a_{\bold{m}}^{g,n}\;
\left(\frac{\partial^3 F_0}{\partial t_0^3}\right)^{1-g-\|\bold{m}\|}\quad\prod_{i=1}^{3g-3+n}
\left(\frac{\partial^{i+3} F_0}{\partial t_0^{i+3}}\right)^{m_i}
$$
with certain constants $a_{\bold{m}}^{g,n}$ depending on $g,n$ and
${\bold{m}}=(m_1,\dots, m_{3g-3+n})$. These constants are linear 
combinations of the intersection numbers $\langle\tau_2^{l_2}\tau_3^{l_3}\dots\rangle$,
and the corresponding $p\,(3g-3+n)\times p\,(3g-3)$ matrix (where $p\,(k)$ is the partition 
number of $k$) can be written quite explicitly (cf. [EYY], Theorem 1 and Appendix D).
Below we give a topological 
interpretation of the numbers $a_{\bold{m}}^{g,n}$; details will appear elsewhere [GoOrZo].

\smallskip
 
Denote by $t_i \in H_{2i}(BU,\bold{Q})$ the class associated with the hyperplane section 
line bundle on $\bold{C}P^i$; the classes $t_i$ form the canonical multiplicative basis of the 
homology ring $H_*(BU,\bold{Q})$. Let us treat the classes $\frac{1}{i!}\,\kappa_{g,n;\,i}\in
H^{2i}(\overline{M}_{g,n},\bold{Q})$ as the components of the Chern character of a (virtual)
vector bundle on $\overline{M}_{g,n}$, and denote by $\alpha_{g,n;\,k}$ the 
corresponding roots.  Put $Q(\alpha)=1+\sum_{i=1}^\infty \alpha^i t_i$ and consider
the homology class in $H_{6g-6+2n}(BU,\bold{Q})$ defined by
$$[\overline{M}_{g,n}]= 
\int_{\overline{M}_{g,n}}\prod_{k}Q(\alpha_{g,n;\,k})=
\sum_{|{\bold{m}}|=3g-3+n}b_{\bold{m}}^{g,n}\,t_1^{m_1}\dots t_{3g-3+n}^{m_{3g-3+n}}.
$$
Similar to $K(x,s_1,\dots)$, we introduce the formal series
$$
B(x,t_1,\dots)=\sum_{g,n}\frac{x^n}{n!}\sum_{|{\bold{m}}|=3g-3+n}b_{\bold{m}}^{g,n}
\prod_{i=1}^{3g-3+n}t_i^{m_i}.
$$

\medskip

{\bf 5.1. Theorem.} {\it (a) The classes $[\overline{M}_{0,n}]$ form
a multiplicative basis of the homology ring $H_*(BU,\bold{Q})$ and
for any $g,n$ we have the decomposition 
$$
[\overline{M}_{g,n}]=\sum_{|{\bold{m}}|=3g-3+n}a_{\bold{m}}^{g,n}
\prod_{i=1}^{3g-3+n}[\overline{M}_{0,n}]^{m_i}
$$
in $H_*(BU,\bold{Q})$ with the same constants $a_{\bold{m}}^{g,n}$ as defined above.

(b) The series $K(x,s_1,\dots)$ and $B(x,t_1,\dots)$ are related 
via the formula
$$
K(x,s_1,\dots)=B(x,t_1,\dots)\left|_{t_i=q_i(s_1,\dots,s_i)}\right. ,
$$
where $q_i$ are the ordinary Schur polynomials defined by
$$
\roman{exp}\,\left(\sum_{i=1}^{\infty}\lambda^is_i\right)=
1+\sum_{j=1}^{\infty}\lambda^j q_j(s_1,\dots ,s_j).
$$
}

\smallskip

The proof essentially follows from Theorem 4.1 and formula (16).

\medskip

{\bf 5.2. Remark.} If we formally put $Q(\psi)=\sum_{i=0}^\infty \psi^i t_i$, 
then for Witten's total free energy we have the formula
$$
F(t_0,t_1,\dots)=\sum_{g,n}\frac{1}{n!}
\int_{\overline{M}_{g,n}}\prod_{i=1}^n Q(\psi_{g,n;i}).
$$
This is reminiscent of J.~Morava's approach to Witten's two dimensional
gravity [Mo], although we do not make use of complex cobordisms.

\medskip

{\bf 6. Asymptotics for the volumes.} We start with defining
some constants. Consider the classical Bessel function
$$
J_0(z)=\sum_{m=0}^{\infty}\frac{(-1)^m}{(m!)^2}\,\left(\frac{z^2}{4}\right)^m.
$$
Denote by $j_0$ its first positive zero and put
$$
x_0=-\frac{1}{2}\,j_0J_0^{\prime}(j_0),\qquad
y_0=-\frac{1}{4}\,j_0^2,\qquad A=-j_0^{-1}J_0^{\prime}(j_0).
\eqno(17)
$$
All  these constants are positive. Put also $V_{g,n}=
\langle \kappa_1^{3g-3+n}\rangle .$ Recall that the
K\"ahler form $\omega_{WP}$ of the Weil--Petersson metric
on $M_{g,n}$ extends as a closed current to $\overline{M}_{g,n}$
and $[\omega_{WP}]=2\pi^2\kappa_1$ in the real cohomology (see [Wo]). 
This means that,
up to a normalization, $V_{g,n}$ are the classical
Weil--Petersson volumes of moduli spaces. In the unstable range
we put by definition $V_{g,n}=0$ for $2g+n\le 2$ and $V_{0,3}=1.$

\medskip

\proclaim{\quad 6.1. Theorem} We have the following asymptotical 
expansions valid for any fixed $g$ and $n\to \infty$:
$$
\frac{V_{g,n}}{n! (n+3g-3)!}\ =\
(n+1)^{\frac{5g-7}{2}}\,x_0^{-n}
\left(B_g + \sum_{k=1}^{\infty}\frac{B_{g,k}}{(n+1)^k}\right) 
\eqno(18)
$$
where
$$
B_0=\frac{1}{A^{\frac{1}{2}}\Gamma\left(-\frac{1}{2}\right)x_0^{\frac{1}{2}}},
\quad B_1=\frac{1}{48},
\eqno(19)
$$
and for $g\ge 2$
$$
B_g=\frac{A^{\frac{g-1}{2}}}{2^{2g-2}\,(3g-3)!\,\Gamma\left(\frac{5g-5}{2}
\right)\,x_0^{\frac{5g-5}{2}}}\;\langle \tau^{3g-3}_2\rangle .
\eqno(20)
$$
\endproclaim
\smallskip

{\bf 6.2. Remark.} By definition,
$$
\langle \tau^{3g-3}_2\rangle =
\int_{\overline{M}_{g,3g-3}}\psi_1^2\dots\psi_{3g-3}^2 .
$$
These numbers can be consecutively calculated using
a simple recursive formula established in [IZu], sec. 6. 
We reproduce it here for a reader's convenience. Put $b_0=-1$,
$b_1= \dfrac{1}{24},$ $b_g=\dfrac{(5g-3)(5g-5)}{2^g\,(3g-3)!}\langle
\tau_2^{3g-3}\rangle $
for $g\ge 2.$ Then we have
$$
b_{g+1}=\frac{25g^2-1}{24}\,b_g +\frac{1}{2}\,
\sum_{m=1}^g b_{g+1-m} b_m \, .
$$
Equivalently, the formal series $\chi (t) =
\sum_{g=0}^{\infty} b_g\,t^{-\frac{5g-1}{2}}$
satisfies the first Painlev\'e equation
$$
\frac{1}{3}\,\chi^{\prime\prime}+\chi^2-t=0.
$$ 

\smallskip

The constants $B_{g,k}$ in the asymptotical expansion (18)
are also calculable in terms of known quantities: see below.

\medskip

{\bf 6.3. Proof of Theorem 6.1.} The proof breaks into two 
logically distinct but tightly interwoven parts.

\smallskip

In the  first part we calculate the smallest
singularity of the function 
$$
\varphi_g(x)=\sum_{n=0}^{\infty} 
\frac{V_{g,n}}{n!\,(n+3g-3)!}\,x^n.
\eqno(21)
$$
and establish its behavior near this singularity.
It turns out that the convergence radius of (21) is $x_0$ (see (17)), and
that near this point the function (or its derivative for $g=0,1$)
can be represented as a  convergent Laurent series in $(x_0-x)^{1/2}$
with the leading term
$$
\varphi_0^{\prime\prime}(x)\sim
-\frac{1}{A^{\frac{1}{2}}}(x_0-x)^{\frac{1}{2}}\, ,
\eqno(22)
$$
$$
\varphi_{1}^{\prime}(x)\sim \frac{1}{48\,(x_0-x)}
\eqno(23)
$$
and  for $g\ge 2$
$$
\varphi_{g}(x) \sim
\frac{A^{\frac{g-1}{2}}}{2^{2g-2} (3g-3)!}\,
\langle\tau_2^{3g-3}\rangle\, (x_0-x)^{-\frac{5g-5}{2}} \,.
\eqno(24)
$$
The right hand side of (18) can now be obtained in the following way.
Expand each monomial  $a_{\nu}(x_0-x)^{\nu}$ into the
Taylor series near $x=0$ and take the sum of the standard asymptotic expansions of their 
coefficients in $n.$
The leading singularities (22)--(24) produce the leading terms in (18).

\smallskip

The justification of this procedure is the second part of the proof.
Instead of proceeding directly, we will show that
we can apply to our situation Theorem 6.1 of [O], Chapter 4,
which is a version of the stationary phase method. 

\smallskip

Now we will explain the precise meaning of (22)--(24)
and prove the statements above.

\smallskip

The genus zero
generating function $y(x):=\varphi_0^{\prime\prime}(x)$
near zero {\it as a formal series} is obtained by inverting
the Bessel function $x(y):=-\sqrt{y}\,J_0^{\prime}(2\sqrt{y}).$
A simple direct proof of this fact based upon
an earlier result by P. Zograf is given in the introduction of [KaMZ].
The derivative of $x(y)$ in $y$ is $J_0(2\sqrt{y}).$
All roots of this function are real and simple (see e.g. [O]).
Clearly, the first root is $y_0$ (see (17)). On the half--line
$(-\infty ,y_0)$ the function $x(y)$ strictly increases
from $-\infty$ to $x_0.$ It follows that these half--lines
are identified by 
a well defined {\it set--theoretic} inverse function $y(x).$
Since $x(y)$ is holomorphic (in fact, entire) and
$y(x)$ is a series in $x$ with positive coefficients,
it converges to the set--theoretic inverse in the interval
$(-x_0,x_0)$ and in the open disk $|x|<x_0.$ Near $y=y_0$
we have $x(y)=x_0-A(y_0-y)^2 + O((y_0-y)^3)$ (see (17)).
Projecting the complex graph of $x(y)$ in $\bold{C}^2$
to the $y$--axis we see that, over a sufficiently small disk
around $y=y_0$, the function $(x_0-x)^{\frac{1}{2}}$ lifted
to this graph is a local parameter at $(x_0,y_0)$. Hence
$y(x)$ can be expanded into Taylor series in $(x_0-x)^{\frac{1}{2}}$
starting with 
$$
y=y_0-\frac{1}{A^{1/2}}\,(x_0-x)^{\frac{1}{2}} +O((x_0-x)).
\eqno(25)
$$
This explains (and proves) (22).

\smallskip

Notice that $x_0$ is the closest to zero singularity only on the
branch of $y(x)$ described above and analytically continued to the
cut complex plane $\bold{C}\setminus [x_0,\infty ].$ The other branches
contain a sequence of real positive ramification points tending to zero.
In fact, the critical values of $J_0(x)$ are real and tend to zero. 

\smallskip

We now pass to the cases $g\ge 1.$
Formula (10) shows that $\varphi_g(x)$ is the specialization of
$K_g(x,s)$ at $s_1=1, s_2=s_3=\dots =0.$ Therefore,
in view of Theorem 4, it is the specialization of
$F_g(t_0,t_1,\dots )$ at  
$t_0=x,\,t_1=0,\, t_k=\dfrac{(-1)^k}{(k-1)!}$ for $k\ge 1.$
In particular, $u_0$ gets replaced by $y(x)$
and prime means now the derivation in $x$ so that
$\varphi_1(x)=\dfrac{1}{24}\,\roman{log}\,y^{\prime}(x)$
which together with (22) establishes (23).

\smallskip

Finally, to treat the case $g\ge 2$ denote by $f_k(x)$
the specialization of $I_k$ corresponding to the above values of
$t_i.$ So, we have
$$
f_2(x)=\frac{y^{\prime\prime}(x)}{y^{\prime}(x)^3},\quad
f_{k+1}=\frac{f_k^{\prime}(x)}{y^{\prime}(x)}\, 
\eqno(26)
$$
and
$$
\varphi_g(x)=\sum_{\sum(i-1)l_i=3g-3} \langle \tau_2^{l_2}\tau_3^{l_3}\dots
\tau_{3g-2}^{l_{3g-2}}\rangle\,
y^{\prime}(x)^{2g-2+\sum l_i}\,\prod_{k=2}^{3g-2}\frac{f_k(x)^{l_k}}{l_k!}\, .
\eqno(27)
$$
Since $y^{\prime}(x)\,x^{\prime}(y)=1$ on the graph of our function,
the r. h. s. of (27) can be rewritten in the following way:
$$
\left. \sum_{\sum(i-1)l_i=3g-3} 
(-1)^{\sum l_i}\langle \tau_2^{l_2}\tau_3^{l_3}\dots
\tau_{3g-2}^{l_{3g-2}}\rangle\,
\frac{1}{x^{\prime}(y)^{2g-2+\sum l_i}}\,
\prod_{k=2}^{3g-2}\frac{x^{(k)}(y)^{l_k}}{l_k!}\,\right|_{y=y(x)} 
\eqno(28)
$$
Each term of (28) as a function of $y$ is meromorphic
on the whole $y$--plane. Hence we can substitute into its
Laurent series at $x_0$
the Taylor series of $y(x)$ at $x_0$ and get the Laurent series
for $\varphi_g(x).$

\smallskip

In particular, an easy induction starting with (25) shows that all $f_k(x)$
are finite at $x_0.$ On the other hand, $y^{\prime}(x)$
replacing $u_0^{\prime}$ in (16) starts with
$\dfrac{1}{2A^{\frac{1}{2}}}\,(x_0-x)^{-\frac{1}{2}}.$
Therefore the leading singularity in (27) is furnished 
by the unique term for which the sum $\sum l_i$ is maximal, that is,
$$
(l_2,\dots ,l_{3g-2})=(3g-3,0,\dots ,0)
$$
To deduce (24) from here, it remains to notice that
$f_2(x_0)=2A.$ 

\smallskip

We now pass to the second part of the proof. Let $q(x)$
denote $y(x)$ or one of the summands in (27).
We want to get the asymptotic expansion (for $n\to\infty $) of
the integral taken over a small circle around zero
$$
\frac{1}{2\pi i}\,\int x^{-(n+1)}\,q(x)dx\, .
\eqno(29)
$$
According to the previous discussion, 
this contour can be continuously deformed in the
definition domain of our branch of $y(x)$.
The natural choice is to take
$\Gamma$ to be the inverse image of the circle $|y|=y_0$
under the map $x\mapsto y(x).$ We owe to Don Zagier
the crucial remark  that $\Gamma$ is a closed
curve without self--intersections, $x\mapsto x(y)$
is a homeomorphism of $\Gamma$, and the circle $|x|=x_0$
lies strictly inside $\Gamma$ except for one point
$x=x_0$ where this circle touches $\Gamma$.
In other words, $|x(y_0e^{i\phi })|>x_0$ unless
$\phi\in 2\pi i\bold{Z}.$ (The easiest way to
convince oneself in this is to look at the
computer generated graph of $|x(y_0e^{i\phi })|-x_0$
produced for us by Don Zagier.
One can also supply a straightforward analytic proof which we omit.)
Denote by $\Gamma_+$ the part of $\Gamma$
in the upper half--plane run over clockwise
and by $\Gamma_-$ the lower part run over counterclockwise
so that
$$
\frac{1}{2\pi i}\,\int x^{-(n+1)}\,q(x)dx\, =
\frac{1}{2\pi i}\,\left(\int_{\Gamma_+} - \int_{\Gamma_-}\right) x^{-(n+1)}\,q(x)dx .
\eqno(30)
$$
Now we have to check that the conditions (I)--(V)
stated in [O], IV.6.1 are satisfied so that we can
apply Theorem 6.1 of loc. cit. to our $\int_{\Gamma_{\pm}}$.
Except for (III), everything is already checked.
In particular, the crucial condition (V) is precisely
the fact that $\Gamma$ includes the circle $|x|=x_0.$
A part of the condition (III) requires the order of $q(x)$
at $x_0$ (Olver's $\lambda -1$) to be greater than $-1$. Since
in our situation this order can be $\dfrac{5}{2}\,(1-g)$
(for the leading term when $g\ge 2$) 
we have to start with representing $q(x)$ as
$q_{\infty}(x)+q_{0}(x)$ where $q_{\infty}$
is a linear combination of negative powers of $(x_0-x)^{\frac{1}{2}}$
and the order of $q_0(x)$ is greater than $-1.$ Then the first
summand is treated directly whereas for the second one
we have the stationary phase expansion.

\smallskip

This completes the proof of Theorem 6.1. We will now complement it
by an explicit prescription for the calculation of all the coefficients
$B_{g,k}$. To adapt Olver's formula IV.(6.19) to our situation,
notice that his $\mu$ is 1, and his $\lambda -1$ is now the order
of $q_0(x)$ at $x=x_0$ which depends on $g$ and on the
summand in (27) that we are treating.  Define the coefficients
$a_k$ by
$$
x_0e^vq_0(x_0e^v)=\sum_{s=0}^{\infty} a_k v^{k+\lambda -1}
\eqno(31)
$$  
near $v=0$ (this is Olver's formula IV.(6.09)). Then the asymptotic expansion 
of (29) for even $g$ reads
$$
\frac{1}{2\pi i}\,\int_{\Gamma_{\pm}} x^{-(n+1)}q_0(x)dx\,= \pm
\frac{1}{2\pi}\,x_0^{-(n+1)}\sum_{k=0}^{\infty}
\Gamma (k+\lambda )\,\frac{a_k}{(n+1)^{k+\lambda}}
\eqno(32)
$$
from which arbitrary  number of coefficients $B_{g,k}$
can be calculated.

\bigskip

{\bf 7. The path groupoid in quantum cohomology.} The first term
of any operad is a groupoid. Since $H_*M$ is
the operad of coalgebras (with comultiplication induced by
the diagonal), we get a bialgebra $B_*=\oplus_{g\ge 1} H_*(\overline{M}_{g,2})$
and the dual bialgebra $B^*=\prod_{g\ge 1} H^*(\overline{M}_{g,2}).$
More precisely, denote here by $b$ the family of boundary
morphisms
$$
b_{g_1,g_2}:\,\overline{M}_{g_1,2}\times \overline{M}_{g_1,2}\to
\overline{M}_{g_1+g_2,2}
$$
glueing the second point of the first curve to the first
point of the second curve.
It induces a multiplication on the homology coalgebra
and a comultiplication on the cohomology algebra
$$
b_*:\,B_*\otimes B_*\to B_*, \quad b^*:\,B^*\to B^*\otimes B^*
$$
making each space a bialgebra.
Moreover, renumbering the structure sections $x_1\leftrightarrow x_2$
we get an involution $s$ of $\overline{M}_{g,2}$ inducing
an involution of $B_*$, resp. $B^*.$ The latter is an automorphism
of the comultiplication on homology, resp. multiplication
on cohomology, but an antiautomorphism of the remaining
two structures, because from the definition one  sees that
$$
b\circ \sigma_{1,2}=s\circ b\circ (s\times s)
$$
where $\sigma_{1,2}$ is the permutation of factors.

\medskip

{\bf 7.1. Proposition.} {\it Let $V$ be a smooth projective
manifold. Then the family of Gromov--Witten correspondences
$I_{g,2}\in A_*(V^2\times \overline{M}_{g,2})$ defines
an action of the algebra $B_*$ and a coaction of
the coalgebra $B^*$ on $H^*(V).$}

\smallskip

This is a part of the general statement that
$H^*(V)$ is a cyclic algebra over the modular operad $H_*M.$

\smallskip

The geometry of this (co)action (explained e.g. in [KoM])
shows that morally it reflects the properties of 
the complexified groupoid of paths on $V$: a stable map of a curve
with two marked points to $V$ should be considered
as a complex path from the first marked point to the second one.
This agrees with the composition by glueing the endpoint of one
path to the starting point of another, and
the smoothing of the resulting singularity plays
the role of homotopy. Notice that there may well
exist non--constant stable maps of genus zero
with only two marked points whereas there are
no stable curves with this property. For this reason,
one may hope that $B_*$ (and possibly $B^*$) 
can be extended by the genus zero component,
perhaps by taking an appropriate (co)homology
of the respective Artin stack $\Cal{M}_{0,2}$
which might be close to the (co)homology of $BG_m.$

\smallskip

As we mentioned in the Introduction, it would
be highly desirable to introduce quantum fundamental
group of $V,$ with appropriate action of our
path groupoid. 
In any case, we are interested in the primitive elements $h$
of $B^*$ satisfying $b^*(h)=h\otimes 1+1\otimes h.$
Clearly, restrictng $L$ to its $n=2$ part we get a supply of such
elements. Their full classification may be easier than that of $L.$ 

\medskip

{\bf Acknowledgements.} We are grateful to Boris Dubrovin for very 
stimulating discussions and Don Zagier for his help
with the proof of Theorem 6.1. PZ thanks MPIM (Bonn) for hospitality
and support during 1998-99 academic year.

\bigskip

\centerline{\bf References}

\medskip

[AC] E.~Arbarello, M.~Cornalba. {\it Combinatorial and
algebro--geometric cohomology classes on the moduli spaces
of curves.} Journ. Alg. Geom. 5 (1996), 705-749.

\smallskip

[EYY] T.~Eguchi, Y.~Yamada, S.--K.~Yang. 
{\it On the genus expansion in the topological string theory.}
Rev. Mod. Phys. 7 (1995), 279.

\smallskip

[FP] C.~Faber, R.~Pandharipande. {\it Hodge integrals
and Gromov--Witten theory.} Preprint math.AG/9810173

\smallskip

[GK] E.~Getzler, M.~M.~Kapranov. {\it Modular operads.}
Comp. Math. 110 (1998), 65-126.

\smallskip

[GoOrZo] V.~Gorbounov, D.~Orlov, P.~Zograf (in preparation).

\smallskip

[IZu] C.~Itzykson, J.--B.~Zuber. {\it Combinatorics
of the modular group II: the Kontsevich integrals.}
Int.~J.~Mod.~Phys. A7 (1992), 5661.

\smallskip

[KabKi] A.~Kabanov, T.~Kimura. {\it Intersection numbers and rank one
cohomological field theories in genus one.} Comm. Math. Phys. 194
(1998), 651-674.

\smallskip

[KaMZ] R.~Kaufmann, Yu.~Manin, D.~Zagier. {\it Higher Weil--Petersson
volumes of moduli spaces of stable $n$--pointed curves.}
Comm. Math. Phys. 181 (1996), 763--787.

\smallskip

[KoM] M.~Kontsevich, Yu.~Manin. {\it Gromov-Witten classes, quantum
cohomology, and enumerative geometry.} Comm. Math. Phys.
164 (1994), 525--562.

\smallskip

[KoMK] M.~Kontsevich, Yu.~Manin (with Appendix by R.~Kaufmann).
{\it Quantum cohomology of a product.} Inv. Math. 124 (1996),
313--340.

\smallskip

[Mo] J.~Morava. {\it Schur Q-functions and a Kontsevich--Witten genus.}
Contemp. Math. 220 (1998), 255--266.

\smallskip

[Mu] D.~Mumford. {\it Towards an enumerative geometry of the moduli
space of curves.} In: Arithmetic and Geometry (M.~Artin and J.~Tate, eds.),
Part II, Birkh\"auser, 1983, 271--328.

\smallskip

[O] F.~W.~J.~Olver. {\it Introduction to asymptotics
and special functions.} Academic Press, 1974.

\smallskip

[W] E.~Witten. {\it Two-dimensional gravity and intersection theory
on moduli space.} Surveys in Diff. Geom. 1 (1991), 243--310.

\smallskip

[Wo] S.~Wolpert. {\it The hyperbolic metric and the geometry of
the universal curve.} J.~Diff.~Geo. 31 (1990), 417--472.

\smallskip

[Zo] P.~Zograf. {\it Weil--Petersson volumes of 
moduli spaces of curves and the genus expansion in
two dimensional gravity.} Preprint math.AG/9811026.

\enddocument